\documentclass[12pt]{article}
\usepackage{amsfonts}
\usepackage{amssymb}
\usepackage{graphicx}
\usepackage{amsmath}

\begin{document}

\begin{center}
{\LARGE Entire Functions Sharing Small Functions With Their Difference
Operators}

\quad

\textbf{Zinel\^{a}abidine} \textbf{LATREUCH}$^{1}$, \textbf{Abdallah}
\textbf{EL FARISSI}$^{2}$ \textbf{and Benharrat} \textbf{BELA\"{I}DI}$^{1}$

\quad

$^{1}$\textbf{Department of Mathematics }

\textbf{Laboratory of Pure and Applied Mathematics }

\textbf{University of Mostaganem (UMAB) }

\textbf{B. P. 227 Mostaganem-(Algeria)}

\textbf{z.latreuch@gmail.com}

\textbf{belaidi@univ-mosta.dz}

$^{2}$\textbf{Department of Mathematics and Informatics, }

\textbf{Faculty of Exact Sciences,}

\textbf{University of Bechar-(Algeria)}

\textbf{elfarissi.abdallah@yahoo.fr}

\quad
\end{center}

\noindent \textbf{Abstract. }We investigate uniqueness problems for an
entire function that shares two small functions of finite order with their
difference operators. In particular, we give a generalization of a result in
$[2]$.

\quad



\section{Introduction and Main Results}

\noindent Throughout this paper, we assume that the reader is familiar with
the fundamental results and the standard notations of the Nevanlinna's value
distribution theory $(\left[ 7\right] ,$ $\left[ 9\right] ,$ $\left[ 12%
\right] )$. In addition, we will use $\rho \left( f\right) $ to denote the
order of growth of $f$ and $\tau \left( f\right) $ to denote the type of
growth of $f$, we say that a meromorphic function $a\left( z\right) $ is a
small function of $f\left( z\right) $ if $T\left( r,a\right) =S\left(
r,f\right) ,$ where $S\left( r,f\right) =o\left( T\left( r,f\right) \right)
, $ as $r\rightarrow \infty $ outside of a possible exceptional set of
finite logarithmic measure, we use $S\left( f\right) $ to denote the family
of all small functions with respect to $f\left( z\right) $. For a
meromorphic function $f\left( z\right) ,$ we define its shift by $%
f_{c}\left( z\right) =f\left( z+c\right) $ $\left( \text{Resp. }f_{0}\left(
z\right) =f\left( z\right) \right) $ and its difference operators by%
\begin{equation*}
\Delta _{c}f\left( z\right) =f\left( z+c\right) -f\left( z\right) ,\text{ \
\ }\Delta _{c}^{n}f\left( z\right) =\Delta _{c}^{n-1}\left( \Delta
_{c}f\left( z\right) \right) ,\text{ }n\in
\mathbb{N}
,\text{ }n\geq 2.
\end{equation*}%
In particular, $\Delta _{c}^{n}f\left( z\right) =\Delta ^{n}f\left( z\right)
$ for the case $c=1.$

\noindent \qquad Let $f\left( z\right) $ and $g\left( z\right) $ be two
meromorphic functions, and let $a\left( z\right) $ be a small function with
respect to $f\left( z\right) $ and $g\left( z\right) .$ We say that $f\left(
z\right) $ and $g\left( z\right) $ share $a\left( z\right) $ CM (counting
multiplicity), provided that $f\left( z\right) -a\left( z\right) $ and $%
g\left( z\right) -a\left( z\right) $ have the same zeros with the same
multiplicities.

\noindent \qquad The problem of meromorphic functions sharing small
functions with their differences is an important topic of uniqueness theory
of meromorphic functions $\left( \text{see, }\left[ 1,4-6\right] \right) $.
In 1986, Jank, Mues and Volkmann $\left( \text{see, }\left[ 8\right] \right)
$ proved:

\quad

\noindent \textbf{Theorem A} \textit{Let }$f$\textit{\ be a nonconstant
meromorphic function, and let }$a\neq 0$\textit{\ be a finite constant. If }$%
f,$\textit{\ }$f^{\prime }$\textit{\ and }$f^{\prime \prime }$\textit{\
share the value }$a$\textit{\ CM, then }$f\equiv f^{\prime }.$

\quad

\noindent In $\left[ 11\right] ,$ P. Li and C. C. Yang gives the following
generalization of Theorem A.

\quad

\noindent \textbf{Theorem B} \textit{Let }$f$\textit{\ be a nonconstant
entire function, let }$a$\textit{\ be a finite nonzero constant, and let }$n$%
\textit{\ be a positive integer. If }$f$\textit{, }$f^{\left( n\right) }$%
\textit{\ and }$f^{\left( n+1\right) }$\textit{\ share the value }$a$\textit{%
\ CM, then }$f\equiv f^{\prime }.$

\textit{\quad }

\noindent \qquad In $\left[ 2\right] ,$ B. Chen et al proved a difference
analogue of result of Theorem A and obtained the following results:

\quad

\noindent \textbf{Theorem C }\textit{Let }$f\left( z\right) $ \textit{be a
nonconstant entire function of finite order, and let }$a\left( z\right)
\left( \not\equiv 0\right) \in S\left( f\right) $\textit{\ be a periodic
entire function with period }$c$\textit{. If }$f\left( z\right) ,$\textit{\ }%
$\Delta _{c}f$\textit{\ and }$\Delta _{c}^{2}f$\textit{\ share }$a\left(
z\right) $\textit{\ CM, then }$\Delta _{c}f\equiv \Delta _{c}^{2}f.$

\quad

\noindent \textbf{Theorem D }\textit{Let }$f\left( z\right) $ \textit{be a
nonconstant entire function of finite order, and let }$a\left( z\right) ,$ $%
b\left( z\right) \left( \not\equiv 0\right) \in S\left( f\right) $\textit{\
be periodic entire functions with period }$c$\textit{. If }$f\left( z\right)
-a\left( z\right) ,$\textit{\ }$\Delta _{c}f\left( z\right) -b\left(
z\right) $\textit{\ and }$\Delta _{c}^{2}f\left( z\right) -b\left( z\right) $%
\textit{\ share }$0$\textit{\ CM, then }$\Delta _{c}f\equiv \Delta
_{c}^{2}f. $

\quad

\noindent \qquad Recently in $\left[ 3\right] ,$ B. Chen and S. Li
generalized Theorem C and proved the following results:

\quad

\noindent \textbf{Theorem E }\textit{Let }$f\left( z\right) $ \textit{be a
nonconstant entire function of finite order, and let }$a\left( z\right)
\left( \not\equiv 0\right) \in S\left( f\right) $\textit{\ be a periodic
entire function with period }$c$\textit{. If }$f\left( z\right) ,$\textit{\ }%
$\Delta _{c}f$\textit{\ and }$\Delta _{c}^{n}f$\textit{\ }$\left( n\geq
2\right) $ \textit{share }$a\left( z\right) $\textit{\ CM, then }$\Delta
_{c}f\equiv \Delta _{c}^{n}f.$

\quad

\noindent \textbf{Theorem F }\textit{Let }$f\left( z\right) $ \textit{be a
nonconstant entire function of finite order. If }$f\left( z\right) ,$\textit{%
\ }$\Delta _{c}f\left( z\right) $\textit{\ and }$\Delta _{c}^{n}f\left(
z\right) $\textit{\ share }$0$\textit{\ CM, then }$\Delta _{c}^{n}f\left(
z\right) =C\Delta _{c}f\left( z\right) ,$\textit{\ where }$C$\textit{\ is a
nonzero constant.}

\quad

\noindent \qquad It is interesting now to see what happening when $f\left(
z\right) $, $\Delta _{c}^{n}f\left( z\right) $\ and $\Delta
_{c}^{n+1}f\left( z\right) $\ $\left( n\geq 1\right) $ share $a\left(
z\right) $ CM. The main of this paper is to give a difference analogue of
result of Theorem B. In fact, we prove that the conclusion of Theorems E and
F remains valid when we replace $\Delta _{c}f\left( z\right) $ by $\Delta
_{c}^{n+1}f\left( z\right) $, and we obtain the following results.

\quad

\noindent \textbf{Theorem 1.1} \textbf{\ }\textit{Let }$f\left( z\right) $
\textit{be a nonconstant entire function of finite order, and let }$a\left(
z\right) \left( \not\equiv 0\right) \in S\left( f\right) $\textit{\ be a
periodic entire function with period }$c$\textit{. If }$f\left( z\right) $%
\textit{, }$\Delta _{c}^{n}f\left( z\right) $\textit{\ and }$\Delta
_{c}^{n+1}f\left( z\right) $\textit{\ }$\left( n\geq 1\right) $ \textit{%
share }$a\left( z\right) $ \textit{CM, then }$\Delta _{c}^{n+1}f\left(
z\right) \equiv \Delta _{c}^{n}f\left( z\right) .$

\quad

\noindent \textbf{Example 1.1 }Let $f\left( z\right) =e^{z\ln 2}$ and $c=1.$
Then, for any $a\in
\mathbb{C}
,$ we notice that $f\left( z\right) ,$ $\Delta _{c}^{n}f\left( z\right) $%
\textit{\ }and $\Delta _{c}^{n+1}f\left( z\right) $\ share $a$\ CM for all $%
n\in
\mathbb{N}
$ and we can easily see that $\Delta _{c}^{n+1}f\left( z\right) \equiv
\Delta _{c}^{n}f\left( z\right) .$ This example satisfies Theorem 1.1.

\quad

\noindent \textbf{Remark 1.1 }In Example 1.1, we have $\Delta
_{c}^{m}f\left( z\right) \equiv \Delta _{c}^{n}f\left( z\right) $ for any
integer $m>n+1.$ However, it remains open when $f\left( z\right) $\textit{, }%
$\Delta _{c}^{n}f\left( z\right) $\textit{\ }and\textit{\ }$\Delta
_{c}^{m}f\left( z\right) $\textit{\ }$\left( m>n+1\right) $ share\textit{\ }$%
a\left( z\right) $ CM, the claim $\Delta _{c}^{n+1}f\left( z\right) \equiv
\Delta _{c}^{n}f\left( z\right) $ in Theorem 1.1 can be replaced by $\Delta
_{c}^{m}f\left( z\right) \equiv \Delta _{c}^{n}f\left( z\right) $ in general.

\quad

\noindent \textbf{Theorem 1.2 }\textit{Let }$f\left( z\right) $\textit{be a
nonconstant entire function of finite order, and let }$a\left( z\right) ,$ $%
b\left( z\right) \left( \not\equiv 0\right) \in S\left( f\right) $\textit{\
be a periodic entire function with period }$c$\textit{. If }$f\left(
z\right) -a\left( z\right) ,$\textit{\ }$\Delta _{c}^{n}f\left( z\right)
-b\left( z\right) $\textit{\ and }$\Delta _{c}^{n+1}f\left( z\right)
-b\left( z\right) $\textit{\ share }$0$\textit{\ CM, then }$\Delta
_{c}^{n+1}f\left( z\right) \equiv \Delta _{c}^{n}f\left( z\right) .$

\quad

\noindent \textbf{Theorem 1.3 }\textit{Let }$f\left( z\right) $\textit{be a
nonconstant entire function of finite order. If }$f\left( z\right) ,$\textit{%
\ }$\Delta _{c}^{n}f\left( z\right) $\textit{\ and }$\Delta
_{c}^{n+1}f\left( z\right) $\textit{\ share }$0$\textit{\ CM, then }$\Delta
_{c}^{n+1}f\left( z\right) \equiv C\Delta _{c}^{n}f\left( z\right) ,$\textit{%
\ where }$C$\textit{\ is a nonzero constant.}

\quad

\noindent \textbf{Example 1.2 }Let $f\left( z\right) =e^{az}$ and $c=1$
where $a\neq 2k\pi i$ $\left( k\in
\mathbb{Z}
\right) ,$ it is clear that $\Delta _{c}^{n}f\left( z\right) =\left(
e^{a}-1\right) ^{n}e^{az}$ for any integer $n\geq 1.$ So, $f\left( z\right)
, $ $\Delta _{c}^{n}f\left( z\right) $\textit{\ }and $\Delta
_{c}^{n+1}f\left( z\right) $\ share $0$\ CM for all $n\in
\mathbb{N}
$ and we can easily see that $\Delta _{c}^{n+1}f\left( z\right) \equiv
C\Delta _{c}^{n}f\left( z\right) $ where $C=e^{a}-1.$ This example satisfies
Theorem 1.3.

\section{Some lemmas}

\noindent \textbf{Lemma 2.1 }$\left[ 10\right] $\ \ \textit{Let }$f$\textit{%
\ and }$g$\textit{\ be meromorphic functions\ such that }$0<$\textit{\ }$%
\rho \left( f\right) ,\rho \left( g\right) <\infty $\textit{\ and }$0<\tau
\left( f\right) ,\tau \left( g\right) <\infty .$\textit{\ Then we have}

\noindent $\left( \text{i}\right) $ \textit{If }$\rho \left( f\right) >\rho
\left( g\right) ,$ \textit{then we obtain}%
\begin{equation*}
\tau \left( f+g\right) =\tau \left( fg\right) =\tau \left( f\right) .
\end{equation*}%
$\left( \text{ii}\right) $ \textit{If }$\rho \left( f\right) =\rho \left(
g\right) $ \textit{and }$\tau \left( f\right) \neq \tau \left( g\right) ,$
\textit{then we get}%
\begin{equation*}
\rho \left( f+g\right) =\rho \left( fg\right) =\rho \left( f\right) =\rho
\left( g\right) .
\end{equation*}%
\textbf{Lemma 2.2 }$\left[ 12\right] $ \textit{Suppose }$f_{j}\left(
z\right) $\textit{\ }$(j=1,2,\cdots ,n+1)$\textit{\ and }$g_{j}\left(
z\right) $\textit{\ }$(j=1,2,\cdots ,n)$\textit{\ }$(n\geq 1)$\textit{\ are
entire functions satisfying the following conditions:}

\noindent $\left( \text{i}\right) $\textit{\ }$\overset{n}{\underset{j=1}{%
\sum }}f_{j}\left( z\right) e^{g_{j}\left( z\right) }\equiv f_{n+1}\left(
z\right) ;$

\noindent $\left( \text{ii}\right) $\textit{\ The order of }$f_{j}\left(
z\right) $\textit{\ is less than the order of }$e^{g_{k}\left( z\right) }$%
\textit{\ for }$1\leq j\leq n+1,$\textit{\ }$1\leq k\leq n.$\textit{\ And
furthermore, the order of }$f_{j}\left( z\right) $\textit{\ is less than the
order of }$e^{g_{h}\left( z\right) -g_{k}\left( z\right) }$\textit{\ for }$%
n\geq 2$\textit{\ and }$1\leq j\leq n+1,$\textit{\ }$1\leq h<k\leq n.$

\noindent \textit{Then }$f_{j}\left( z\right) \equiv 0,$\textit{\ }$\left(
j=1,2,\cdots n+1\right) .$

\quad

\noindent \textbf{Lemma 2.3 }$\left[ 5\right] $ \textit{Let }$c\in
\mathbb{C}
,$\textit{\ }$n\in
\mathbb{N}
,$\textit{\ and let }$f\left( z\right) $\textit{\ be a meromorphic function
of finite order. Then for any small periodic function }$a\left( z\right) $%
\textit{\ with period }$c,$\textit{\ with respect to }$f\left( z\right) ,$%
\begin{equation*}
m\left( r,\frac{\Delta _{c}^{n}f}{f-a}\right) =S\left( r,f\right) ,
\end{equation*}%
\textit{where the exceptional set associated with }$S\left( r,f\right) $%
\textit{\ is of at most finite logarithmic measure.}

\section{Proof of the Theorems}

\noindent \textbf{Proof of the Theorem 1.1.} Suppose on the contrary to the
assertion that $\Delta _{c}^{n}f\left( z\right) \not\equiv \Delta
_{c}^{n+1}f\left( z\right) .$ Note that $f\left( z\right) $ is a nonconstant
entire function of finite order. By Lemma 2.3, for $n\geq 1$, we have%
\begin{equation*}
T\left( r,\Delta _{c}^{n}f\right) =m\left( r,\Delta _{c}^{n}f\right) \leq
m\left( r,\frac{\Delta _{c}^{n}f}{f}\right) +m\left( r,f\right) \leq T\left(
r,f\right) +S\left( r,f\right) .
\end{equation*}%
Since $f\left( z\right) $, $\Delta ^{n}f\left( z\right) $\ and $\Delta
^{n+1}f\left( z\right) $\ $\left( n\geq 1\right) $ share $a\left( z\right) $
CM, then%
\begin{equation}
\frac{\Delta _{c}^{n}f\left( z\right) -a\left( z\right) }{f\left( z\right)
-a\left( z\right) }=e^{P\left( z\right) }  \tag{3.1}
\end{equation}%
and
\begin{equation}
\frac{\Delta _{c}^{n+1}f\left( z\right) -a\left( z\right) }{f\left( z\right)
-a\left( z\right) }=e^{Q\left( z\right) },  \tag{3.2}
\end{equation}%
where $P$ and $Q$ are polynomials. Set%
\begin{equation}
\varphi \left( z\right) =\frac{\Delta _{c}^{n+1}f\left( z\right) -\Delta
_{c}^{n}f\left( z\right) }{f\left( z\right) -a\left( z\right) }.  \tag{3.3}
\end{equation}%
From $\left( 3.1\right) $ and $\left( 3.2\right) ,$ we get $\varphi \left(
z\right) =e^{Q\left( z\right) }-e^{P\left( z\right) }.$ Then, by supposition
and $(3.3)$, we see that $\varphi \left( z\right) \not\equiv 0$. By Lemma
2.3, we deduce that
\begin{equation}
T\left( r,\varphi \right) =m\left( r,\varphi \right) \leq m\left( r,\frac{%
\Delta _{c}^{n+1}f}{f-a}\right) +m\left( r,\frac{\Delta _{c}^{n}f}{f-a}%
\right) +O\left( 1\right) =S\left( r,f\right) .  \tag{3.4}
\end{equation}%
Note that $\frac{e^{Q\left( z\right) }}{\varphi \left( z\right) }-\frac{%
e^{P\left( z\right) }}{\varphi \left( z\right) }=1.$ By using the second
main theorem and $(3.4)$, we have%
\begin{equation*}
T\left( r,\frac{e^{Q}}{\varphi }\right) \leq \overline{N}\left( r,\frac{e^{Q}%
}{\varphi }\right) +\overline{N}\left( r,\frac{\varphi }{e^{Q}}\right) +%
\overline{N}\left( r,\frac{1}{\frac{e^{Q}}{\varphi }-1}\right) +S\left( r,%
\frac{e^{Q}}{\varphi }\right)
\end{equation*}%
\begin{equation*}
=\overline{N}\left( r,\frac{e^{Q}}{\varphi }\right) +\overline{N}\left( r,%
\frac{\varphi }{e^{Q}}\right) +\overline{N}\left( r,\frac{\varphi }{e^{P}}%
\right) +S\left( r,\frac{e^{Q}}{\varphi }\right)
\end{equation*}%
\begin{equation}
=S\left( r,f\right) +S\left( r,\frac{e^{Q}}{\varphi }\right) .  \tag{3.5}
\end{equation}%
Thus, by $(3.4)$ and $(3.5)$, we have $T(r,e^{Q})$ $=S(r,f)$. Similarly, $%
T(r,e^{P})=S(r,f)$. Setting now $g\left( z\right) =f\left( z\right) -a\left(
z\right) ,$ we have from $\left( 3.1\right) $ and $\left( 3.2\right) $%
\begin{equation}
\Delta _{c}^{n}g\left( z\right) =g\left( z\right) e^{P\left( z\right)
}+a\left( z\right)  \tag{3.6}
\end{equation}%
and%
\begin{equation}
\Delta _{c}^{n+1}g\left( z\right) =g\left( z\right) e^{Q\left( z\right)
}+a\left( z\right) .  \tag{3.7}
\end{equation}%
By $\left( 3.6\right) $ and $\left( 3.7\right) ,$ we have%
\begin{equation*}
g\left( z\right) e^{Q\left( z\right) }+a\left( z\right) =\Delta _{c}\left(
\Delta _{c}^{n}g\left( z\right) \right) =\Delta _{c}\left( g\left( z\right)
e^{P\left( z\right) }+a\left( z\right) \right) .
\end{equation*}%
Thus%
\begin{equation*}
g\left( z\right) e^{Q\left( z\right) }+a\left( z\right) =g_{c}\left(
z\right) e^{P_{c}\left( z\right) }-g\left( z\right) e^{P\left( z\right) },
\end{equation*}%
which implies%
\begin{equation}
g_{c}\left( z\right) =M\left( z\right) g\left( z\right) +N\left( z\right) ,
\tag{3.8}
\end{equation}%
where $M\left( z\right) =e^{-P_{c}\left( z\right) }\left( e^{P\left(
z\right) }+e^{Q\left( z\right) }\right) $ and $N\left( z\right) =a\left(
z\right) e^{-P_{c}\left( z\right) }.$ From $\left( 3.8\right) ,$ we have
\begin{equation*}
g_{2c}\left( z\right) =M_{c}\left( z\right) g_{c}\left( z\right)
+N_{c}\left( z\right) =M_{c}\left( z\right) \left( M\left( z\right) g\left(
z\right) +N\left( z\right) \right) +N_{c}\left( z\right) ,
\end{equation*}%
hence
\begin{equation*}
g_{2c}\left( z\right) =M_{c}\left( z\right) M_{0}\left( z\right) g\left(
z\right) +N^{1}\left( z\right) ,
\end{equation*}%
where $N^{1}\left( z\right) =M_{c}\left( z\right) N_{0}\left( z\right)
+N_{c}\left( z\right) .$ By the same method, we can deduce that%
\begin{equation}
g_{ic}\left( z\right) =\left( \underset{k=0}{\overset{i-1}{\prod }}%
M_{kc}\left( z\right) \right) g\left( z\right) +N^{i-1}\left( z\right) \text{
}\left( i\geq 1\right) ,  \tag{3.9}
\end{equation}%
where $N^{i-1}\left( z\right) $ $\left( i\geq 1\right) $ is an entire
function depending on $a\left( z\right) ,e^{P\left( z\right) },e^{Q\left(
z\right) }$ and their differences. Now, we can rewrite $\left( 3.6\right) $
as%
\begin{equation}
\overset{n}{\underset{i=1}{\sum }}C_{n}^{i}\left( -1\right)
^{n-i}g_{ic}\left( z\right) =\left( e^{P\left( z\right) }-\left( -1\right)
^{n}\right) g\left( z\right) +a\left( z\right) .  \tag{3.10}
\end{equation}%
By $\left( 3.9\right) $ and $\left( 3.10\right) ,$ we have%
\begin{equation*}
\overset{n}{\underset{i=1}{\sum }}C_{n}^{i}\left( -1\right) ^{n-i}\left(
\left( \underset{k=0}{\overset{i-1}{\prod }}M_{kc}\left( z\right) \right)
g\left( z\right) +N^{i-1}\left( z\right) \right) -\left( e^{P\left( z\right)
}-\left( -1\right) ^{n}\right) g\left( z\right) =a\left( z\right)
\end{equation*}%
which implies%
\begin{equation}
A\left( z\right) g\left( z\right) +B\left( z\right) =0,  \tag{3.11}
\end{equation}%
where
\begin{equation*}
A\left( z\right) =\overset{n}{\underset{i=1}{\sum }}C_{n}^{i}\left(
-1\right) ^{n-i}\underset{k=0}{\overset{i-1}{\prod }}M_{kc}\left( z\right)
-e^{P\left( z\right) }+\left( -1\right) ^{n}
\end{equation*}%
and
\begin{equation*}
B\left( z\right) =\overset{n}{\underset{i=1}{\sum }}C_{n}^{i}\left(
-1\right) ^{n-i}N^{i-1}\left( z\right) -a\left( z\right) .
\end{equation*}%
It is clear that $A\left( z\right) $ and $B\left( z\right) $ are small
functions with respect to $f\left( z\right) .$ If\textbf{\ }$A\left(
z\right) \not\equiv 0$, then $\left( 3.11\right) $ yields the contradiction
\begin{equation*}
T\left( r,f\right) =T\left( r,g\right) =T\left( r,\frac{B}{A}\right)
=S\left( r,f\right) .
\end{equation*}%
Suppose now that $A\left( z\right) \equiv 0,$ rewrite the equation $A\left(
z\right) \equiv 0$ as%
\begin{equation*}
\overset{n}{\underset{i=1}{\sum }}C_{n}^{i}\left( -1\right) ^{n-i}\underset{%
k=0}{\overset{i-1}{\prod }}e^{-P_{\left( k+1\right) c}}\left(
e^{P_{kc}}+e^{Q_{kc}}\right) =e^{P}-\left( -1\right) ^{n}.
\end{equation*}%
We can rewrite the left side of above equality as%
\begin{equation*}
\overset{n}{\underset{i=1}{\sum }}C_{n}^{i}\left( -1\right) ^{n-i}e^{-%
\overset{i}{\underset{k=1}{\sum }}P_{kc}}\underset{k=0}{\overset{i-1}{\prod }%
}\left( e^{P_{kc}}+e^{Q_{kc}}\right)
\end{equation*}%
\begin{equation*}
=\overset{n}{\underset{i=1}{\sum }}C_{n}^{i}\left( -1\right) ^{n-i}e^{-%
\overset{i}{\underset{k=1}{\sum }}P_{kc}}e^{\overset{i-1}{\underset{k=0}{%
\sum }}P_{kc}}\underset{k=0}{\overset{i-1}{\prod }}\left(
1+e^{Q_{kc}-P_{kc}}\right)
\end{equation*}%
\begin{equation*}
\overset{n}{=\underset{i=1}{\sum }}C_{n}^{i}\left( -1\right)
^{n-i}e^{P-P_{ic}}\underset{k=0}{\overset{i-1}{\prod }}\left(
1+e^{Q_{kc}-P_{kc}}\right) .
\end{equation*}%
So%
\begin{equation}
\overset{n}{\underset{i=1}{\sum }}C_{n}^{i}\left( -1\right)
^{n-i}e^{P-P_{ic}}\underset{k=0}{\overset{i-1}{\prod }}\left(
1+e^{h_{kc}}\right) =e^{P}-\left( -1\right) ^{n},  \tag{3.12}
\end{equation}%
where $h_{kc}=Q_{kc}-P_{kc}.$ On the other hand, let $\Omega _{i}=\left\{
0,1,\cdots ,i-1\right\} $ be a finite set of $i$ elements, and
\begin{equation*}
P\left( \Omega _{i}\right) =\{\varnothing ,\left\{ 0\right\} ,\left\{
1\right\} ,\cdots ,\left\{ i-1\right\} ,\left\{ 0,1\right\} ,\left\{
0,2\right\} ,\cdots ,\Omega _{i}\},
\end{equation*}%
where $\varnothing $ is an empty set. It is easy to see that%
\begin{equation*}
\underset{k=0}{\overset{i-1}{\prod }}\left( 1+e^{h_{kc}}\right) =1+\underset{%
A\in P\left( \Omega _{i}\right) \backslash \left\{ \varnothing \right\} }{%
\sum }\exp \left( \underset{j\in A}{\sum }h_{jc}\right)
\end{equation*}%
\begin{equation}
=1+\left[ e^{h}+e^{h_{c}}+\cdots +e^{h_{\left( i-1\right) c}}\right] +\left[
e^{h+h_{c}}+e^{h+h_{2c}}+\cdots \right] +\cdots +\left[ e^{h+h_{c}+\cdots
+h_{\left( i-1\right) c}}\right] .  \tag{3.13}
\end{equation}%
Dividing the proof on two parts:

\noindent \textbf{Part (1). }$h\left( z\right) $ is non-constant polynomial.
Suppose that $h\left( z\right) =a_{m}z^{m}+\cdots +a_{0}$ $\left( a_{m}\neq
0\right) ,$ since $P\left( \Omega _{i}\right) \subset P\left( \Omega
_{i+1}\right) ,$ then by $\left( 3.12\right) $ and $\left( 3.13\right) $ we
have%
\begin{equation*}
\overset{n}{\underset{i=1}{\sum }}C_{n}^{i}\left( -1\right)
^{n-i}e^{P-P_{ic}}+\alpha _{1}e^{a_{m}z^{m}}+\alpha
_{2}e^{2a_{m}z^{m}}+\cdots +\alpha _{n}e^{na_{m}z^{m}}=e^{P}-\left(
-1\right) ^{n}
\end{equation*}%
which is equivalent to%
\begin{equation}
\alpha _{0}+\alpha _{1}e^{a_{m}z^{m}}+\alpha _{2}e^{2a_{m}z^{m}}+\cdots
+\alpha _{n}e^{na_{m}z^{m}}=e^{P},  \tag{3.14}
\end{equation}%
where $\alpha _{i}$ $\left( i=0,\cdots ,n\right) $ are entire functions of
order less than $m.$ Moreover,%
\begin{equation*}
\alpha _{0}=\overset{n}{\underset{i=1}{\sum }}C_{n}^{i}\left( -1\right)
^{n-i}e^{P-P_{ic}}+\left( -1\right) ^{n}
\end{equation*}%
\begin{equation*}
=e^{P}\left( \overset{n}{\underset{i=1}{\sum }}C_{n}^{i}\left( -1\right)
^{n-i}e^{-P_{ic}}+\left( -1\right) ^{n}e^{-P}\right) =e^{P}\Delta
_{c}^{n}e^{-P}.
\end{equation*}%
$\left( \text{i}\right) $ If $\deg P>m$, then we obtain from $\left(
3.14\right) $ that
\begin{equation*}
\deg P\leq m
\end{equation*}%
which is a contradiction.

\noindent $\left( \text{ii}\right) $ If $\deg P<m,$ then by using Lemma 2.1
and $\left( 3.14\right) $ we obtain
\begin{equation*}
\deg P=\rho \left( e^{P}\right) =\rho \left( \alpha _{0}+\alpha
_{1}e^{a_{m}z^{m}}+\alpha _{2}e^{2a_{m}z^{m}}+\cdots +\alpha
_{n}e^{na_{m}z^{m}}\right) =m,
\end{equation*}%
which is also a contradiction.

\noindent $\left( \text{iii}\right) $ If $\deg P=m,$ then we suppose that $%
P\left( z\right) =dz^{m}+P^{\ast }\left( z\right) $ where $\deg P^{\ast }<m.$
We have to study two subcases:

\noindent $\left( \ast \right) $ If $d\neq ia_{m}$ $\left( i=1,\cdots
,n\right) ,$ then we have
\begin{equation*}
\alpha _{1}e^{a_{m}z^{m}}+\alpha _{2}e^{2a_{m}z^{m}}+\cdots +\alpha
_{n}e^{na_{m}z^{m}}-e^{P^{\ast }}e^{dz^{m}}=-\alpha _{0}.
\end{equation*}%
By using Lemma 2.2, we obtain $e^{P^{\ast }}\equiv 0,$ which is impossible.

\noindent $\left( \ast \ast \right) $ Suppose now that there exists at most $%
j\in \left\{ 1,2,\cdots ,n\right\} $ such that $d=ja_{m}.$ Without loss of
generality, we assume that $j=n.$ Then we rewrite $\left( 3.14\right) $ as%
\begin{equation*}
\alpha _{1}e^{a_{m}z^{m}}+\alpha _{2}e^{2a_{m}z^{m}}+\cdots +\left( \alpha
_{n}-e^{P^{\ast }}\right) e^{na_{m}z^{m}}=-\alpha _{0}.
\end{equation*}%
By using Lemma 2.2, we have $\alpha _{0}\equiv 0,$ so $\Delta
_{c}^{n}e^{-P}=0.$ Thus%
\begin{equation}
\overset{n}{\underset{i=0}{\sum }}C_{n}^{i}\left( -1\right)
^{n-i}e^{-P_{ic}}\equiv 0.  \tag{3.15}
\end{equation}%
Suppose that $\deg P=\deg h=m>1$ and
\begin{equation*}
P\left( z\right) =b_{m}z^{m}+b_{m-1}z^{m-1}+...+b_{0},\text{ }\left(
b_{m}\neq 0\right) .
\end{equation*}
Note that for $j=0,1,\cdots ,n,$ we have
\begin{equation*}
P\left( z+jc\right) =b_{m}z^{m}+\left( b_{m-1}+mb_{m}jc\right) z^{m-1}+\beta
_{j}\left( z\right) ,
\end{equation*}%
where $\beta _{j}\left( z\right) $ are polynomials with degree less than $%
m-1.$ Rewrite $\left( 3.15\right) $ as%
\begin{equation*}
e^{-\beta _{n}\left( z\right) }e^{-b_{m}z^{m}-\left( b_{m-1}+mb_{m}nc\right)
z^{m-1}}-ne^{-\beta _{n-1}\left( z\right) }e^{-b_{m}z^{m}-\left(
b_{m-1}+mb_{m}\left( n-1\right) c\right) z^{m-1}}
\end{equation*}%
\begin{equation}
+\cdots +\left( -1\right) ^{n}e^{-\beta _{0}\left( z\right)
}e^{-b_{m}z^{m}-b_{m-1}z^{m-1}}\equiv 0.  \tag{3.16}
\end{equation}%
For any $0\leq l<k\leq n,$ we have
\begin{equation*}
\rho \left( e^{-b_{m}z^{m}-\left( b_{m-1}+mb_{m}lc\right) z^{m-1}-\left(
-b_{m}z^{m}-\left( b_{m-1}+mb_{m}kc\right) z^{m-1}\right) }\right) =\rho
\left( e^{-mb_{m}\left( l-k\right) cz^{m-1}}\right)
\end{equation*}%
\begin{equation*}
=m-1,
\end{equation*}%
and for $j=0,1,\cdots ,n,$ we see that
\begin{equation*}
\rho \left( e^{\beta _{j}}\right) \leq m-2.
\end{equation*}%
By this, together with $\left( 3.16\right) $ and Lemma 2.2, we obtain $%
e^{-\beta _{n}\left( z\right) }\equiv 0,$ which is impossible. Suppose now
that $P\left( z\right) =\mu z+\eta $ $\left( \mu \neq 0\right) $ and $%
Q\left( z\right) =\alpha z+\beta $ because if $\deg Q>1,$ then we back to
the case $\left( \text{ii}\right) .$ It easy to see that%
\begin{equation*}
\Delta _{c}^{n}e^{-P}=\overset{n}{\underset{i=0}{\sum }}C_{n}^{i}\left(
-1\right) ^{n-i}e^{-\mu \left( z+ic\right) -\eta }=e^{-P}\overset{n}{%
\underset{i=0}{\sum }}C_{n}^{i}\left( -1\right) ^{n-i}e^{-\mu ic}
\end{equation*}%
\begin{equation*}
=e^{-P}\left( e^{-\mu c}-1\right) ^{n}.
\end{equation*}%
This together with $\Delta _{c}^{n}e^{-P}\equiv 0$ gives $\left( e^{-\mu
c}-1\right) ^{n}\equiv 0,$ which yields $e^{\mu c}\equiv 1.$ Therefore, for
any $j\in
\mathbb{Z}
$%
\begin{equation*}
e^{P\left( z+jc\right) }=e^{\mu z+\mu jc+\eta }=\left( e^{\mu c}\right)
^{j}e^{P\left( z\right) }=e^{P\left( z\right) }.
\end{equation*}%
In order to prove that $e^{Q\left( z\right) }$ is also periodic entire
function with period $c,$ we suppose the contrary, which means that $%
e^{\alpha c}\neq 1$. Since $e^{P\left( z\right) }$ is of period $c,$ then by
$\left( 3.14\right) $, we get%
\begin{equation}
\alpha _{1}e^{\left( \alpha -\mu \right) z}+\alpha _{2}e^{2\left( \alpha
-\mu \right) z}+\cdots +\alpha _{n}e^{n\left( \alpha -\mu \right) z}=e^{\mu
z+\eta },  \tag{3.17}
\end{equation}%
where $\alpha _{i}$ $\left( i=1,\cdots ,n\right) $ are constants. In
particular,%
\begin{equation*}
\alpha _{n}=e^{n\left( \beta -\eta \right) +\alpha c\frac{n\left( n-1\right)
}{2}}
\end{equation*}%
and
\begin{equation*}
\alpha _{1}=\left[ \underset{i=1}{\overset{n}{\sum }}C_{n}^{i}\left(
-1\right) ^{n-i}+\underset{i=2}{\overset{n}{\sum }}C_{n}^{i}\left( -1\right)
^{n-i}e^{\alpha c}\right.
\end{equation*}%
\begin{equation*}
\left. +\underset{i=3}{\overset{n}{\sum }}C_{n}^{i}\left( -1\right)
^{n-i}e^{2\alpha c}+\cdots +e^{\left( n-1\right) \alpha c}\right] e^{\left(
\beta -\eta \right) }
\end{equation*}%
\begin{equation*}
=[C_{n}^{1}\left( -1\right) ^{n-1}+C_{n}^{2}\left( -1\right) ^{n-2}\left(
1+e^{\alpha c}\right) +C_{n}^{3}\left( -1\right) ^{n-3}\left( 1+e^{\alpha
c}+e^{2\alpha c}\right)
\end{equation*}%
\begin{equation*}
+\cdots +C_{n}^{n}\left( -1\right) ^{n-n}\left( 1+e^{\alpha c}+\cdots
+e^{\left( n-1\right) \alpha c}\right) ]e^{\left( \beta -\eta \right) }
\end{equation*}%
\begin{equation*}
=[C_{n}^{1}\left( -1\right) ^{n-1}\frac{e^{\alpha c}-1}{e^{\alpha c}-1}%
+C_{n}^{2}\left( -1\right) ^{n-2}\frac{e^{2\alpha c}-1}{e^{\alpha c}-1}%
+C_{n}^{3}\left( -1\right) ^{n-3}\frac{e^{3\alpha c}-1}{e^{\alpha c}-1}
\end{equation*}%
\begin{equation*}
+\cdots +C_{n}^{n}\left( -1\right) ^{n-n}\frac{e^{n\alpha c}-1}{e^{\alpha
c}-1}]e^{\left( \beta -\eta \right) }
\end{equation*}%
\begin{equation*}
=[C_{n}^{1}\left( -1\right) ^{n-1}\left( e^{\alpha c}-1\right)
+C_{n}^{2}\left( -1\right) ^{n-2}\left( e^{2\alpha c}-1\right)
+C_{n}^{3}\left( -1\right) ^{n-3}\left( e^{3\alpha c}-1\right)
\end{equation*}%
\begin{equation*}
+\cdots +C_{n}^{n}\left( -1\right) ^{n-n}\left( e^{n\alpha c}-1\right) ]%
\frac{e^{\left( \beta -\eta \right) }}{e^{\alpha c}-1}
\end{equation*}%
\begin{equation*}
=\left[ \underset{i=0}{\overset{n}{\sum }}C_{n}^{i}\left( -1\right)
^{n-i}e^{i\alpha c}-\left( -1\right) ^{n}-\underset{i=1}{\overset{n}{\sum }}%
C_{n}^{i}\left( -1\right) ^{n-i}\right] \frac{e^{\left( \beta -\eta \right) }%
}{e^{\alpha c}-1}
\end{equation*}%
\begin{equation*}
=\left( e^{\alpha c}-1\right) ^{n-1}e^{\left( \beta -\eta \right) }.
\end{equation*}%
Rewrite $\left( 3.17\right) $ as
\begin{equation}
\alpha _{1}e^{\left( \alpha -2\mu \right) z}+\alpha _{2}e^{\left( 2\alpha
-3\mu \right) z}+\cdots +\alpha _{n}e^{\left( n\alpha -\left( n+1\right) \mu
\right) z}=e^{\eta },  \tag{3.18}
\end{equation}%
it is clear that for each $1\leq l<m\leq n,$ we have%
\begin{equation*}
\rho \left( e^{\left( m\alpha -\left( m+1\right) \mu -l\alpha +\left(
l+1\right) \mu \right) z}\right) =\rho \left( e^{\left( m-l\right) \left(
\alpha -\mu \right) z}\right) =1.
\end{equation*}%
We have the following two cases:

\noindent $\left( \text{i}_{1}\right) $ If $j\alpha -\left( j+1\right) \mu
\neq 0$ for all $j\in \left\{ 1,2,\cdots ,n\right\} ,$ which means that
\begin{equation*}
\rho \left( e^{\left( j\alpha -\left( j+1\right) \mu \right) z}\right) =1,%
\text{ }1\leq j\leq n
\end{equation*}%
then, by applying Lemma 2.2 we obtain $e^{\eta }\equiv 0,$ which is a
contradiction.

\noindent $\left( \text{i}_{2}\right) $ If there exists $\left( \text{at
most one}\right) $ an integer $j\in \left\{ 1,2,\cdots ,n\right\} $ such
that $j\alpha -\left( j+1\right) \mu =0.$ Without loss of generality, assume
that $e^{\left( n\alpha -\left( n+1\right) \mu \right) z}=1,$ the equation $%
\left( 3.18\right) $ will be
\begin{equation*}
\alpha _{1}e^{\left( \alpha -2\mu \right) z}+\alpha _{2}e^{\left( 2\alpha
-3\mu \right) z}+\cdots +\alpha _{n-1}e^{\left( \left( n-1\right) \alpha
-n\mu \right) z}=e^{\eta }-e^{n\left( \beta -\eta \right) +\alpha c\frac{%
n\left( n-1\right) }{2}}
\end{equation*}%
and by applying Lemma 2.2, we obtain $\alpha _{1}=\left( e^{\alpha
c}-1\right) ^{n-1}e^{\left( \beta -\eta \right) }\equiv 0,$ which is
impossible. So, by $\left( \text{i}_{1}\right) $ and $\left( \text{i}%
_{2}\right) ,$ we deduce that $e^{\alpha c}\equiv 1$. Therefore, for any $%
j\in
\mathbb{Z}
$ we have
\begin{equation*}
e^{Q\left( z+jc\right) }=e^{\alpha z+\beta }\left( e^{\alpha c}\right)
^{j}=e^{Q\left( z\right) },
\end{equation*}%
which implies that $e^{Q}$ is periodic of period $c.$ Since $e^{P\left(
z\right) }$ is of period $c,$ then by $\left( 3.1\right) ,$ we obtain
\begin{equation}
\Delta _{c}^{n+1}f\left( z\right) =e^{P}\Delta _{c}f\left( z\right) ,
\tag{3.19}
\end{equation}%
then $\Delta _{c}^{n+1}f\left( z\right) $ and $\Delta _{c}f\left( z\right) $
share $0$ CM. Substituting $\left( 3.19\right) $ into the second equation $%
\left( 3.2\right) ,$ we get%
\begin{equation}
e^{P\left( z\right) }\Delta _{c}f\left( z\right) =e^{Q\left( z\right)
}\left( f\left( z\right) -a\left( z\right) \right) +a\left( z\right) .
\tag{3.20}
\end{equation}%
Since $e^{P\left( z\right) }$ and $e^{Q\left( z\right) }$ are of period $c,$
then by $\left( 3.20\right) ,$ we obtain%
\begin{equation}
\Delta _{c}^{n+1}f\left( z\right) =e^{Q-P}\Delta _{c}^{n}f\left( z\right) .
\tag{3.21}
\end{equation}%
So, $\Delta ^{n+1}f\left( z\right) $ and $\Delta ^{n}f\left( z\right) $
share $0,a\left( z\right) $ CM, combining $\left( 3.1\right) ,$ $\left(
3.2\right) $ and $\left( 3.21\right) ,$ we deduce that
\begin{equation*}
\frac{\Delta ^{n+1}f\left( z\right) -a\left( z\right) }{\Delta ^{n}f\left(
z\right) -a\left( z\right) }=\frac{\Delta ^{n+1}f\left( z\right) }{\Delta
^{n}f\left( z\right) },
\end{equation*}%
and we get
\begin{equation*}
\Delta ^{n+1}f\left( z\right) =\Delta ^{n}f\left( z\right)
\end{equation*}%
which is a contradiction. Suppose now that $P=c_{1}$ and $Q=c_{2}$ are
constants $\left( e^{c_{1}}\neq e^{c_{2}}\right) .$ By $\left( 3.8\right) $
we have%
\begin{equation*}
g_{c}\left( z\right) =\left( e^{c_{2}-c_{1}}+1\right) g\left( z\right)
+a\left( z\right) e^{-c_{1}}
\end{equation*}%
by the same%
\begin{equation*}
g_{2c}\left( z\right) =\left( e^{c_{2}-c_{1}}+1\right) ^{2}g\left( z\right)
+a\left( z\right) e^{-c_{1}}\left( \left( e^{c_{2}-c_{1}}+1\right) +1\right)
.
\end{equation*}%
By induction, we obtain%
\begin{equation*}
g_{nc}\left( z\right) =\left( e^{c_{2}-c_{1}}+1\right) ^{n}g\left( z\right)
+a\left( z\right) e^{-c_{1}}\underset{i=0}{\overset{n-1}{\sum }}\left(
e^{c_{2}-c_{1}}+1\right) ^{i}
\end{equation*}%
\begin{equation*}
=\left( e^{c_{2}-c_{1}}+1\right) ^{n}g\left( z\right) +a\left( z\right)
e^{-c_{2}}\left( \left( e^{c_{2}-c_{1}}+1\right) ^{n}-1\right) .
\end{equation*}%
Rewrite the equation $\left( 3.6\right) $ as%
\begin{equation*}
\Delta _{c}^{n}g\left( z\right) =\overset{n}{\underset{i=0}{\sum }}%
C_{n}^{i}\left( -1\right) ^{n-i}\left[ \left( e^{c_{2}-c_{1}}+1\right)
^{i}g\left( z\right) +a\left( z\right) e^{-c_{2}}\left( \left(
e^{c_{2}-c_{1}}+1\right) ^{i}-1\right) \right]
\end{equation*}%
\begin{equation*}
=e^{c_{1}}g\left( z\right) +a\left( z\right) .
\end{equation*}%
Since $A\left( z\right) \equiv 0,$ then we have
\begin{equation*}
\overset{n}{\underset{i=0}{\sum }}C_{n}^{i}\left( -1\right) ^{n-i}\left(
e^{c_{2}-c_{1}}+1\right) ^{i}=e^{c_{1}}
\end{equation*}%
and%
\begin{equation*}
\overset{n}{\underset{i=0}{\sum }}C_{n}^{i}\left( -1\right) ^{n-i}\left(
\left( e^{c_{2}-c_{1}}+1\right) ^{i}-1\right) =e^{c_{2}}
\end{equation*}%
which are equivalent to%
\begin{equation*}
e^{n\left( c_{2}-c_{1}\right) }=e^{c_{1}}
\end{equation*}%
and%
\begin{equation*}
e^{n\left( c_{2}-c_{1}\right) }=e^{c_{2}}
\end{equation*}%
which is a contradiction.

\noindent \textbf{Part (2). }$h\left( z\right) $ is a constant. We show
first that $P\left( z\right) $ is a constant. If $\deg P>0,$ from the
equation $\left( 3.12\right) ,$ we see
\begin{equation*}
\deg P\leq \deg P-1,
\end{equation*}%
which is a contradiction. Then $P\left( z\right) $ must be a constant and
since $h\left( z\right) =Q\left( z\right) -P\left( z\right) $ is a constant,
we deduce that both of $P\left( z\right) $ and $Q\left( z\right) $ is
constant. This case is impossible too (the last case in Part (1)), and we
deduced that $h\left( z\right) $ can not be a constant. Thus, the proof of
Theorem 1.1 is completed.

\quad

\noindent \textbf{Proof of the Theorem 1.2. }Setting $g\left( z\right)
=f\left( z\right) +b\left( z\right) -a\left( z\right) ,$ we can remark that
\begin{equation*}
g\left( z\right) -b\left( z\right) =f\left( z\right) -a\left( z\right) ,
\end{equation*}%
\begin{equation*}
\Delta _{c}^{n}g\left( z\right) -b\left( z\right) =\Delta _{c}^{n}f\left(
z\right) -b\left( z\right)
\end{equation*}%
and
\begin{equation*}
\Delta _{c}^{n+1}g\left( z\right) -b\left( z\right) =\Delta _{c}^{n}f\left(
z\right) -b\left( z\right) ,\text{ }n\geq 2.
\end{equation*}%
Since $f\left( z\right) -a\left( z\right) ,$\textit{\ }$\Delta
_{c}^{n}f\left( z\right) -b\left( z\right) $\textit{\ }and $\Delta
_{c}^{n+1}f\left( z\right) -b\left( z\right) $\ share $0$\ CM, it follows
that $g\left( z\right) ,$\textit{\ }$\Delta _{c}^{n}g\left( z\right) $%
\textit{\ }and $\Delta _{c}^{n+1}g\left( z\right) $ share $b\left( z\right) $
CM. By using Theorem 1.1, we deduce that $\Delta _{c}^{n+1}g\left( z\right)
\equiv \Delta _{c}^{n}g\left( z\right) ,$ which leads to $\Delta
_{c}^{n+1}f\left( z\right) \equiv \Delta _{c}^{n}f\left( z\right) $ and the
proof of Theorem 1.2 is completed.

\quad

\noindent \textbf{Proof of the Theorem 1.3. }Note that $f\left( z\right) $
is a nonconstant entire function of finite order. Since $f\left( z\right) ,$%
\textit{\ }$\Delta _{c}^{n}f\left( z\right) $\textit{\ }and $\Delta
_{c}^{n+1}f\left( z\right) $\ share $0$\textit{\ }CM, then we have%
\begin{equation}
\frac{\Delta _{c}^{n}f\left( z\right) }{f\left( z\right) }=e^{P\left(
z\right) }  \tag{3.22}
\end{equation}%
and
\begin{equation}
\frac{\Delta _{c}^{n+1}f\left( z\right) }{f\left( z\right) }=e^{Q\left(
z\right) },  \tag{3.23}
\end{equation}%
where $P$ and $Q$ are polynomials. If $Q-P$ is a constant, then we can get
easily from $\left( 3.22\right) $ and $\left( 3.23\right) $%
\begin{equation*}
\Delta _{c}^{n+1}f\left( z\right) =e^{Q\left( z\right) -P\left( z\right)
}\Delta _{c}^{n}f\left( z\right) :\equiv C\Delta _{c}^{n}f\left( z\right) .
\end{equation*}%
This complete our proof. If $Q-P$ is a not constant, with a similar arguing
as in the proof of Theorem 1.1, we can deduce that the case $\deg P=\deg
\left( Q-P\right) >1$ is impossible. For the case $\deg P=\deg \left(
Q-P\right) =1,$ we can obtain that $e^{P\left( z\right) }$ is periodic
entire function with period $c.$ This together with $\left( 3.22\right) $
yields%
\begin{equation}
\Delta _{c}^{n+1}f\left( z\right) =e^{P\left( z\right) }\Delta _{c}f\left(
z\right)  \tag{3.24}
\end{equation}%
which means that $f\left( z\right) ,$ $\Delta _{c}f\left( z\right) $ and $%
\Delta _{c}^{n+1}f\left( z\right) $ share $0$ CM. Thus, by Theorem F, we
obtain%
\begin{equation*}
\Delta _{c}^{n+1}f\left( z\right) \equiv C\Delta _{c}f\left( z\right)
\end{equation*}%
which is a contradiction with $\left( 3.22\right) $ and $\deg P=1.$ Theorem
1.3 is thus proved.

\quad

\noindent \textbf{Acknowledgements.} The authors are grateful to the referee
for his/her valuable comments which lead to the improvement of this paper.

\begin{center}
{\Large References}
\end{center}

\noindent $\left[ 1\right] \ $W. Bergweiler, J. K. Langley, \textit{Zeros of
differences of meromorphic functions}, Math. Proc. Cambridge Philos. Soc.
142 (2007), no. 1, 133--147.

\noindent $\left[ 2\right] \ $B. Chen, Z. X. Chen and S. Li,\textit{\
Uniqueness theorems on entire functions and their difference operators or
shifts}, Abstr. Appl. Anal. 2012, Art. ID 906893, 8 pp.

\noindent $\left[ 3\right] \ $B. Chen, and S. Li, \textit{Uniquness problems
on entire functions that share a small function with their difference
operators}, Adv. Difference Equ. 2014, 2014:311, 11 pp.

\noindent $\left[ 4\right] \ $Y. M. Chiang, S. J. Feng, \textit{On the
Nevanlinna characteristic of }$f\left( z+\eta \right) $ \textit{and
difference equations in the complex plane, }Ramanujan J. 16 (2008), no. 1,
105-129.

\noindent $\left[ 5\right] \ $R. G. Halburd, R. J. Korhonen, \textit{%
Difference analogue of the lemma on the logarithmic derivative with
applications to difference equations, }J. Math. Anal. Appl. 314 (2006)%
\textit{, }no. 2, 477-487.

\noindent $\left[ 6\right] \ $R. G. Halburd, R. J. Korhonen, \textit{%
Nevanlinna theory for the difference operator}, Ann. Acad. Sci. Fenn. Math.
31 (2006), no. 2, 463--478.

\noindent $\left[ 7\right] \ $W. K. Hayman, \textit{Meromorphic functions},
Oxford Mathematical Monographs Clarendon Press, Oxford 1964.

\noindent $\left[ 8\right] \ $G. Jank. E. Mues and L. Volkmann, \textit{%
Meromorphe Funktionen, die mit ihrer ersten und zweiten Ableitung einen
endlichen Wert teilen}, Complex Variables Theory Appl. 6 (1986), no. 1,
51--71.

\noindent $\left[ 9\right] \ $I. Laine, \textit{Nevanlinna theory and
complex differential equations}, de Gruyter Studies in Mathematics, 15.
Walter de Gruyter \& Co., Berlin, 1993.

\noindent $\lbrack 10]$ Z. Latreuch and B. Bela\"{\i}di, \textit{Estimations
about the order of growth and the type of meromorphic functions in the
complex plane}, An. Univ. Oradea, Fasc. Matematica, Tom XX (2013), Issue No.
1, 179-186.

\noindent $\lbrack 11]$ P. Li and C. C. Yang, \textit{Uniqueness theorems on
entire functions and their derivatives}, J. Math. Anal. Appl. 253 (2001),
no. 1, 50--57.

\noindent $\left[ 12\right] \ $C. C. Yang, H. X. Yi, \textit{Uniqueness
theory of meromorphic functions}, Mathematics and its Applications, 557.
Kluwer Academic Publishers Group, Dordrecht, 2003.

\end{document}